\documentclass[11pt]{article}

\usepackage{amsmath,amssymb,lineno}

\newtheorem{theorem}{Theorem}
\newtheorem{ass}{Theorem}[section]

\newtheorem{lemma}[ass]{Lemma}

\newtheorem{definition}[ass]{Definition}

\newcommand{\qed}{\hspace*{\fill} \rule{7pt}{7pt}}
\newcommand{\Proof}{\noindent{\bf Proof.}\ \ }
\newcommand{\T}{\mathcal T}
\newcommand{\G}{\mathcal G}
\newcommand{\C}{\mathcal C}

\begin{document}

\title{The effect of local majority on global majority in connected graphs}

\author{
Yair Caro
\thanks{Department of Mathematics, University of Haifa at Oranim, Tivon
36006, Israel. Email: yacaro@kvgeva.org.il}
\and
Raphael Yuster
\thanks{Department of Mathematics, University of Haifa, Haifa
31905, Israel. Email: raphy@math.haifa.ac.il}
}

\date{}

\maketitle

\setcounter{page}{1}

\begin{abstract}
Let $\G$ be an infinite family of connected graphs and let $k$ be a positive integer.
We say that $k$ is {\em forcing} for $\G$ if for all $G \in \G$ but finitely many, the following holds.
Any $\{-1,1\}$-weighing of the edges of $G$ for which all connected subgraphs on $k$ edges are positively weighted implies that $G$ is positively weighted.
Otherwise, we say that it is {\em weakly forcing} for $\G$ if any such weighing implies that the weight of $G$ is bounded from below by a constant.
Otherwise we say that $k$ {\em collapses} for $\G$.

We classify $k$ for some of the most prominent classes of graphs,
such as all connected graphs, all connected graphs with a given maximum degree and all connected graphs with a given average degree.
\end{abstract}

\section{Introduction}

Voting majority problems are an important research topic in combinatorics and computer science.
In the graph-theoretic setting of these problems, a binary voting function (i.e a function
into $\{-1,1\}$) is assigned to either the vertices or the edges of the graph according
to some set of rules. The meta-question is what can then be deduced about the overall outcome of the vote in the entire graph?
See \cite{ACFGPO-2015,berger-2001,BHHM-b1995,FHL-1986,FM-1999,peleg-2002,PZ-2010,woodall-1992,xu-2001} for a representative
sample of some very natural problems in this area.

In this paper we consider a voting problem on the edges of a connected graph according to the
very natural rule that all connected subgraphs of a bounded given size vote positively.
In particular, our results generalize some of the results in \cite{CHM-2016} who
considered a voting function on a path. Problems showing the effect (or lack thereof) of properties of all small connected subgraphs on
the same property for the whole graph appear throughout combinatorics. Perhaps the best known example is the result of Erd\H{o}s \cite{erdos-1959} showing that
a graph with relatively large chromatic number may have all small connected subgraphs acyclic.

For a graph $G$ and a weighing $w:E(G) \rightarrow \{-1,1\}$, let $w(G)=\sum_{e \in E(G)} w(e)$ denote the total weight of $G$.
We say that $G$ is {\em positively weighted} (or just {\em positive}) if $w(G) > 0$. The weight of a subgraph is the sum of the weights of its edges.
Clearly, for every $k \ge 1$, if every subgraph on $k$ edges of $G$ is positive, then $G$ is positive as well.

Suppose that all {\em connected} $k$-edge subgraphs of $G$ are positive. Can we then deduce that $G$ is positive?
This is trivially true for $k=1$ but already for $k=2$ one needs to impose a connectivity requirement on $G$ in order to avoid trivial obstructions.
Indeed, if $G$ is just a matching, then it has no connected $k$-edge subgraph for $k \ge 2$ so the aforementioned condition of every connected $k$-edge subgraph being positive,
which vacuously holds, does not imply that $G$ is positive. Hence, hereafter, we shall assume that $G$ is connected.
Once connectivity is assumed, any weighing $w$ with the property that all connected $2$-edge subgraphs are positive, implies that $G$ is indeed positive,
unless $G$ is $K_2$.
It is also not difficult to prove (see the next section) that any weighing $w$ with the property that all connected $3$-edge subgraphs are positive, implies that $G$ is non-negatively weighted, unless $G \in \{K_2,K_{1,2}\}$. 
The same holds for $k=4$ unless $G \in \{K_2,K_{1,2},K_3,P_4,K_{1,3}\}$ since if every connected $4$-edge subgraph is positive,
then so is every connected $3$-edge subgraph. In fact, for $k=4$, $G$ will be positive.
What happens, then, for larger $k$? This turns out to exhibit a rather interesting, and sometimes dichotomous behavior, depending on the class of graphs in question.

To be a bit more formal we introduce some notation. We say that a weighing $w$ as above is {\em $k$-local positive} if $w$ is positive on every connected $k$-edge subgraph.
So our question can be rephrased as asking for which graphs (or, more qualitatively, for which families of graphs) is it true that $k$-local positivity implies positivity.
If positivity is not implied, we can ask how negative could $w(G)$ be when $w$ is $k$-local positive. Will it just be bounded by some number independent of the
size of the graphs in the family, or would it be unbounded? This suggests the following definitions.

\begin{definition}[forcing, weakly forcing, collapse]
Let $\G$ be an infinite family of connected graphs and let $k$ be a positive integer.
We say that $k$ is {\em forcing} for $\G$ if for all $G \in \G$ but finitely many, the following holds:
Any $k$-local positive weighing of $G$ implies that $G$ is positive.\\
If $G$ is not forcing, we say that $k$ is {\em weakly forcing} for $\G$ if for all $G \in \G$ but finitely many, the following holds:
Any $k$-local positive weighing $w$ of $G$ has $w(G) \ge C$ where $C$ is a (possibly negative) constant.
In case $k$ is weakly forcing, we may seek the maximum possible such $C$ which we denote by $f_{\G}(k)$ (and observe that $f_{\G}(k) \le 0$).\\
If $k$ is not weakly forcing, we say that $k$ {\em collapses} for $\G$. In this case we define $f_{\G}(k) = -\infty$.
Here we may seek the supremum over all $c \le 0$ such that $w(G) \ge c|E(G)|-o(|E(G)|)$ for all $G \in \G$. We denote this supremum by $c_{\G}(k)$.
\end{definition}

We use the following notation for some of the most popular families of connected graphs.
For a family of graphs $\G$, let $\G_\Delta$ denote its subset consisting of all elements with maximum degree at most $\Delta$
and let $\G_{\overline \alpha}$ denote its subset consisting of all elements with average degree at least $\alpha$.
Let $\C$ be the family of all connected graphs and let $\T$ be the family of all trees.
Our first main result is a classification of the above for every $k$ and for every family $\G$ with $\C \supseteq \G \supseteq \T$.
\begin{theorem}
\label{t:1}
Let $\C \supseteq \G \supseteq \T$. Then
\begin{enumerate}
\item
$k=1,2,4$ are forcing for $\G$.
\item
$k=3,5,6,8$ are weakly forcing for $\G$ with $f_\G(3)=0$, $f_\G(5)=-1$, $f_\G(6)=0$, $f_G(8)=-1$.
\item
$k=7$ and all $k \ge 9$ collapse for $G$. Furthermore:
\begin{enumerate}
\item
If $k = 4m+3$, then 											 $c_\G(k)  = -\frac{m}{3m+2}$\;.
\item
If $k = 4m+2$, then 											 $c_\G(k)  = -\frac{m-1}{3m+1}$\;.
\item
If $k = 4m+1$, then $-\frac{m-1}{3m+1} \ge c_\G(k) \ge -\frac{2m-1}{6m+1}$\;.
\item
If $k = 4m$,   then $-\frac{m-2}{3m}   \ge c_\G(k) \ge -\frac{2m-3}{6m-1}$\;.
\end{enumerate}
\end{enumerate}
\end{theorem}
Our second main result is a classification of the above for every $k$, every $\Delta$, and for every family $\G_\Delta$ with $\C_\Delta \supseteq \G_\Delta \supseteq \T_\Delta$.
\begin{theorem}
\label{t:2}
Let $\Delta$ be a positive integer and let $\C_\Delta \supseteq \G_\Delta \supseteq \T_\Delta$.
Every $k$ is forcing for $\G_2$. For every $\Delta \ge 3$ we have:
\begin{enumerate}
\item
$k=1,2,3,4,6$ are forcing for $\G_\Delta$.
\item
$k=5,8$ are weakly forcing for $\G_\Delta$ with $f_{\G_\Delta}(5)=-1$, $f_{\G_\Delta}(8)=-1$.
\item
$k=7$ and all $k \ge 9$ collapse for $G_\Delta$. Furthermore:
\begin{enumerate}
\item
If $k = 4m+3$, then $-\frac{2m-1}{6m+5} \ge c_{\G_\Delta}(k) \ge -\frac{m}{3m+2}$\;.
\item
If $k = 4m+2$, then $-\frac{2m-3}{6m+3} \ge c_{\G_\Delta}(k) \ge -\frac{m-1}{3m+1}$\;.
\item
If $k = 4m+1$, then $-\frac{m-1}{3m+1} \ge c_{\G_\Delta}(k) \ge -\frac{2m-1}{6m+1}$\;.
\item
If $k = 4m$,   then $-\frac{m-2}{3m}   \ge c_{\G_\Delta}(k) \ge -\frac{2m-3}{6m-1}$\;.
\end{enumerate}
\end{enumerate}
\end{theorem}
Finally, we consider connected graphs with average degree at least $\alpha$.
It is not difficult to prove that for $\alpha$ sufficiently large, all $k$ are forcing for $\C_{\overline \alpha}$.
We determine the exact threshold when this happens.
\begin{theorem}\label{t:3}
For all $\alpha < 8$, infinitely many $k$ collapse for $\C_{\overline \alpha}$.
For all $\alpha \ge 8$, all $k$ are forcing for $\C_{\overline \alpha}$.
\end{theorem}

The rest of this paper is organized as follows. We first prove in Section 2 that in many cases, classifying $k$ for a family of graphs reduces to classifying $k$ for
its sub-family of trees. In the three sections following it, we prove our main results, Theorems \ref{t:1},\ref{t:2},\ref{t:3}, respectively.
The final section consists of some concluding notes and open problems.

\section{On families that contain all trees}

We start with the proof of the following property satisfied by all families of graphs which are between $\C$ and $\T$ (families of connected graphs that contain all trees).
\begin{theorem}\label{t:trees-graphs}
Let $k \ge 1$ be an integer. The following hold for all $\T \subseteq \G \subseteq \C$.
\begin{enumerate}
\item
$k$ is forcing for $\T$ if and only if $k$ is forcing for  $\G$. Furthermore, if $k$ is forcing for $\G$, then for every $t > 0$, and  for
all $G \in  \G$ but finitely many, it holds that any $k$-local positive weighing of $G$ satisfies $w(G) \ge t$.
\item
$k$ is weakly forcing for $\T$ if and only if $k$ is weakly forcing for $\G$ and $f_{\T}(k)=f_{\G}(k)$.
\item
$k$ collapses for $\T$ if and only if $k$ collapses for $\G$ and $c_{\T}(k)=c_{\G}(k)$
\end{enumerate}
The same statements hold for all $\T_\Delta \subseteq \G_\Delta \subseteq \C_\Delta$ for every $\Delta \ge k/2-1$. 
\end{theorem}
\Proof
Observe that since $\T \subseteq \G$, we have that if $k$ is forcing for $\G$, then $k$ is forcing for $\T$.
We prove the converse. Suppose $k$ is forcing for $\T$. Then, for all  $T \in \T$ but finitely many, it holds that $k$-local positivity implies $w(T) > 0$.
Hence, for some $k_0=k_0(k)$ it holds that all trees with at least $k_0$ edges satisfy $w(T) > 0$.
We will show that all graphs in $G \in \G$ with at least $k_0+1$ vertices (i.e. all but finitely many), a $k$-local positive weighing implies $w(G) > 0$.

We apply induction on the number of edges of $G$ where the base case is $k_0$. In this case, $G$ must be a tree since it is connected and has at least $k_0+1$ vertices (thus precisely $k_0+1$ vertices)
and the statement holds by assumption.
Assume it holds for graphs with $n-1 \ge k_0$ edges and we prove that it holds for graphs with $n$ edges. If some edge with positive weight is on a cycle, then we can remove it and remain with a
connected graph with $n-1$ edges, which satisfies $k$-local positivity, so by induction we are done.

Otherwise, all cycles consist only of negative edges. In other words, if we repeatedly replace each cycle of length $s$ with a star
$K_{1,s}$ where the $s$ leaves play the role
of the cycle vertices and the center vertex is new (this procedure is well known as the  $Y-\Delta$ transform, see \cite{CIM-1998,truemper-1989}), we will obtain a tree with $n$ edges,
where the total weight is unchanged, since at each step we replaced a totally negative cycle of length $s$ with a totally negative $K_{1,s}$.
We claim that the new tree satisfies $k$-local positivity. In fact, we show that after each step of a $Y-\Delta$ transform, $k$-local positivity remains intact. Consider $G$ before a certain
$Y-\Delta$ transform, and suppose we replace some negative cycle $C_s$ of $G$ with a negative $K_{1,s}$ and obtained a graph $G'$.
Notice that if $G'$ has a  non-positive connected $k$-edge subgraph $H$, then we can assume that at least one edge of the new $K_{1,s}$ participates in $H$ since otherwise $H$ would already have been
a connected subgraph of $G$. But now we claim that all edges of the $K_{1,s}$ participate in $H$. Indeed, suppose some edge $e$ does not participate. Then, we can remove from $H$ an edge that is not
on the $K_{1,s}$ and such that after removal, the obtained $(k-1)$-edge subgraph remains connected, and then add $e$ to remain with a connected $k$-edge subgraph and the weight remains negative (it may decrease or stay the same). But now, if all edges of the $K_{1,s}$ are used in $H$, then in $G$ we can replace the edges of the $K_{1,s}$ with the edges of the cycle $C_s$ and obtain
a connected $k$-edge subgraph with non-positive weight in $G$, a contradiction to our assumption.
Having proved that the tree satisfies $k$-local positivity we have that $w(G)=w(T) > 0$, as required.

Notice further that in our $Y-\Delta$ transform, original vertices never have their degree increased, but the center of the new $K_{1,s}$ is a new vertex of degree $s$. But also notice that
$G$ cannot have negative cycles of length larger than $\lceil k/2-1 \rceil$ as otherwise it violates $k$-local positivity. Therefore, if $\Delta \ge k/2-1$, such a transform on a graph in $\G_\Delta$ 
results also in a graph in $\G_\Delta$. Hence the same consequence holds for $\G_\Delta$ as stated in the theorem.

Assume now that $k$ is forcing for $\G$. We prove that for every $t > 0$, and for all but finitely many elements of $\G$ it holds that any
$k$-local positive weighing of $G \in \G$ satisfies $w(G) \ge t$. As shown above, it suffices to prove the claim for the family $\T$ of all trees. So, again let $k_0$ be the least integer such that for all trees with at least $k_0$ edges,
every $k$-local positive weighing of $T$ implies $w(T) > 0$. It is a folklore well-known result that any tree with at least two edges can be partitioned
into two sub-trees, each one containing at least $n/3$ edges. So, a tree with $3k_0$ edges has at least two edge-disjoint sub-trees, each of which containing at least $k_0$ edges.
So, a $k$-local positive weighing of such a tree induces a $k$-local positive weighing of each part, hence $w(T) \ge 2$.
Inductively on $t$, we obtain that if the tree has $3(t-1)k_0$ edges or more, it has $w(T) \ge t$. In fact, as we see, $w(T)$ grows linearly with the number of edges of $T$.

For the second part of the theorem, notice first that if $k$ is weakly forcing for $G$ then it is trivially weakly forcing for $\T$ and $f_\G(k) \le f_\T(k)$.
The converse, namely if $k$ is weakly forcing for $\T$, then it is weakly forcing for $\G$ and furthermore $f_\T(k) \le f_\G(k)$ is proved by exactly the
same inductive argument as in the first part. So does the case for $\G_\Delta$ versus $\T_\Delta$.

The third part of the theorem immediately follows from the first two parts, or directly using the same inductive proof.
\qed

\section{Classifying $k$ for all  trees}

Let $\G$ denote a family of connected graphs that contains all trees (i.e $\T \subseteq \G$). The previous section shows that in order to classify which $k$ are forcing, weakly forcing, or collapse
for $\G$, we only need to classify $k$ for $\T$. We begin with the following lemma.
\begin{lemma}\label{l:1234}
Each of $k=1,2,4$ is forcing for $\T$.
On the other hand $k=3$ is weakly forcing for $\T$ with $f_{\T}(3)=0$. Nevertheless, for every $\Delta$, $k=3$ is forcing for $\T_\Delta$.
\end{lemma}
\Proof
The cases $k=1,2$ are trivial.
We prove that $k=4$ is forcing for $\T$.
Consider some weighing of a tree $T$ which is $4$-local positive and assume that $T$ has at least $4$ edges.
Notice that the set of edges receiving weight $-1$ forms an induced matching,
as otherwise there would have been a connected $4$-edge subgraph having non-positive weight.
Furthermore, each edge with weight $-1$ has the property that all its adjacent edges have weight $1$ and since each
negative edge has at least one positive edge adjacent only to it but not to any other  negative edge, we have that the number of
edges with weight $1$ is not smaller than the number of edges with weight $-1$, and that equality can only be obtained if each of the edges with weight $-1$ has one endpoint
as a leaf and the other endpoint of degree $2$. It follows that equality can only be obtained when $G=S_{t,2}$ is obtained from a star $S_t$ (the star with $t \ge 2$ edges) by subdividing each edge
and the weighing assigns $-1$ to the $t$ leaf edges and $1$ to the $t$ non-leaf edges. But such a weighing has a connected $4$-edge subgraph (in fact, a $4$-edge path) which is
not positive.

For the case $k=3$, notice that the aforementioned proof for $k=4$ carries through assuming $G$ has at least three edges and the only obstruction are the $S_{t,2}$ and the weighing that assigns
$-1$ to the $t$ leaf edges and $1$ to the $t$ non-leaf edges. But even in this case, the overall weight is zero,
showing that $k=3$ is weakly forcing for $\T$ and $f_{\T}(3)=0$. Furthermore, in $\T_\Delta$ the only obstructions of this form are $S_{t,2}$ for $t \le \Delta$, hence
only finitely many obstructions. Thus, by definition, $k=3$ is forcing for $\T_\Delta$.
\qed

The cases $k=5,6,8$ are somewhat more involved, as can be seen from the proof of the next two lemmata.
\begin{lemma}\label{l:k56}
Each of $k=5,6$ is weakly forcing for $\T$ with $f_{\T}(5)=-1$, $f_{\T}(6)=0$.
\end{lemma}
\Proof
We begin with constructions of infinite families of trees showing that $f_{\T}(5)\le -1$ and $f_{\T}(6) \le 0$.
Consider the star $S_t$ and subdivide each edge $r-1$ times (so each edge is replaced with a path on $r$ edges). Denote the resulting tree by $S_{t,r}$.
Consider $S_{t,4}$ and assign weight $-1$ to the leaf edges and to the edges adjacent to them, and assign $1$ to the other edges. This forms a $6$-local positive weighing but the overall
weight is $0$.
Now remove one leaf edge from $S_{t,4}$. This yields one branch on three edges while keeping the other $t-1$ branches intact with four edges.
Call the new tree $R_{t,4}$. Again assign weight $-1$ to the leaf edges and to the edges adjacent to them, and assign $1$ to the other edges.
This forms a $5$-local positive weighing of $R_{t,4}$ but the overall weight is $-1$.

We next prove that $f_{\T}(5) \ge -1$ by showing that for every tree $T$ with at least $5$ edges and for a given $5$-local positive weighing of $T$, it holds that $w(T) \ge -1$.
The proof is by induction on the number of edges where the base case are all trees with $5 \le n \le 8$ edges.
We first verify that the base case indeed holds.
The cases $n=5,6,7$ are trivial since we cannot have less than three positive edges.
The case $n=8$ holds since we cannot have less than four positive edges.
To see this, suppose we had only three positive edges (clearly we have at least three) and remove one of them from $T$. If $T$ remains connected, we are done as for the case $n=7$.
If it is disconnected, one of the components has at least three negative edges and hence $T$ contains a connected $5$-edge subgraph with negative weight.

Now, suppose $T$ is a tree with $n \ge 9$ edges. 
We may assume that all leaf edges of $T$ are negative, as otherwise we can use the same inductive proof as in the previous paragraph.
We may further assume that $T$ has no three adjacent negative edges as otherwise there would be a connected $5$-edge subgraph with negative weight, violating our assumption.
Let $r$ be an arbitrary non-leaf vertex of $T$. We root $T$ at $r$ obtaining a rooted tree. We will use the standard terminology of rooted trees: parent, sibling, child, level.
Vertex $r$ is at level $0$ while a vertex $v$ is at level $i$ if its parent is in level $i-1$. Let $v_1$ be a leaf of largest possible level $t$
and let $v_{i+1}$ be the parent of $v_i$ for $i=1,\ldots,t$ where $v_{t+1}=r$.
Notice that $v_1$ has at most one sibling, since if it had more, we would have three leaf edges incident to $v_2$, and thus three adjacent negative edges, violating
our assumption.

We therefore consider two cases where the first case is that $v_1$ has a unique sibling $v_0$. In this case we must have $t \ge 2$ so $v_3$ exists.
Now, if $(v_2,v_3)$ is negative we would have three negative adjacent edges, violating our assumption. Hence $(v_2,v_3)$ is positive.
In fact, this now implies that all edges incident with $v_3$ are positive since otherwise if $e$ is a negative edge incident with $v_3$, then $e$ together with $(v_1,v_2)$, $(v_0,v_2)$, $(v_2,v_3)$
would have total weight $-2$ and we would have a connected $5$-edge subgraph with weight at most $-1$.
If $v_3$ has no children other than $v_2$, then $t \ge 3$ so $v_4$ exists. But then removing $v_0,v_1,v_2,v_3$ and hence edges $(v_0,v_2),(v_1,v_2),(v_2,v_3),(v_3,v_4)$ results
in a tree with $n-4$ edges whose weight is the same as the weight of $T$, so we are done by induction.
Now, if $v_3$ has children other than $v_2$, let $u$ be such a child.
But then $(v_3,u)$ has positive weight so $u$ is not a leaf. Let $x$ be a child of $u$. Then $x$ must be a leaf since its level is the same as the level of $v_1$ which was maximum.
But since all leaf edges are negative, $(u,x)$ is negative. But now $(x,u),(u,v_3),(v_2,v_3),(v_1,v_2)(v_0,v_2)$ is a connected $5$-edge subgraph with negative weight, violating our assumption.

We remain with the case that $v_1$ has no siblings. Now, if $(v_2,v_3)$ is positive, then removing vertices $v_1,v_2$ and hence edges
$(v_1,v_2)$ and $(v_2,v_3)$ results in a tree with $n-2$ edges whose weight is the same as the weight of $T$, so we are done by induction.
We may therefore assume that $(v_2,v_3)$ is negative. So all other edges incident with $v_3$ must be positive. Thus, as in the previous case, $v_3$ cannot have children other than $v_2$.
Hence $v_4$ exists and $(v_3,v_4)$ is positive. Once again, all other edges incident with $v_4$ must be positive. Furthermore, $v_4$ cannot have children other than $v_3$.
So, $v_5$ exists and $(v_4,v_5)$ is positive. But now we can delete vertices $v_1,v_2,v_3,v_4$ and hence edges $(v_1,v_2),(v_2,v_3),(v_3,v_4),(v_4,v_5)$ resulting
in a tree with $n-4$ edges whose weight is the same as the weight of $T$, so we are done by induction.

Finally, observe that whenever we have used induction, we have either removed a single edge (whose weight is positive), two edges (one positive and one negative) or four edges
(two positive and two negative), and since we have proved our base case for the four consecutive values $5,6,7,8$, we are done.

The proof for $k=6$ showing that $f_T(6)=0$ is the same in its inductive part. We just need to take care of the base case of the induction
which will now be all trees on $n=6,7,8,9$ edges. We verify that the base cases indeed holds, namely that $6$-local positive weighing implies nonnegative weighing.
The cases $n=6,7,8$ are trivial since we cannot have less than four positive edges.
The case $n=9$ holds since we cannot have less than five positive edges.
To see this, suppose we had only four positive edges (clearly we have at least four) and remove one of them from $T$. If $T$ remains connected, we are done as for the case $n=8$.
If it is disconnected, one of the components has at least three negative edges and hence $T$ contains a connected $6$-edge subgraph with non-positive weight.

\begin{lemma}\label{l:k8}
$k=8$ is weakly forcing for $\T$ with $f_{\T}(8)=-1$.
\end{lemma}
\Proof
Let $R_{t,6}$ be obtained from $S_{t,6}$ by removing a leaf edge. Assign weight $-1$ to the leaf edges, to the edges adjacent to them, and to the edges adjacent to the latter.
Assign $1$ to the other edges. This forms an $8$-local positive weighing of $R_{t,6}$ but the overall weight is $-1$.

We prove that $f_{\T}(8) \ge -1$ by showing that for every tree $T$ with at least $8$ edges and for a given $8$-local positive weighing of $T$, it holds that $w(T) \ge -1$.
In this case, we need to establish the induction base for six cases, $n=8,9,10,11,12,13$. We verify that the base cases hold, namely that $8$-local positive weighing implies total weight at most $-1$.
Observe that this is trivial for $n=8,9,10,11$ as we must have at least five positive edges. For $n=12$ we claim that we need at least six positive edges. Indeed, otherwise by removing
a positive edge we have a connected component with at least four negative edges and at most four positive edges, hence there is a connected $8$-edge subgraph which is not positive.
For $n=13$ we also need at least six positive edges for the exact same reason.

As for the inductive step, we use the same notations as before. As in the previous lemma, we may assume that all leaf edges are negative, that no four negative edges are pairwise adjacent, 
that we have rooted the tree at some vertex $r$ where $v_1$ is a leaf of lowest possible level $t$, and the path to the root is $v_1,v_2,\ldots,v_{t+1}=r$.

Here there are three cases according to the number of siblings of $v_1$ which may be $0,1,2$. We prove the inductive step for the case where $v_1$ has no siblings. The proofs of the other cases are similar (and are, in fact, easier as we have a ``head start'' of more than one negative edge). So, $(v_1,v_2)$ is negative, and if $(v_2,v_3)$ is positive we are done by removing these two edges and using induction for $n-2$.

So, from this point onwards we can assume that $(v_1,v_2)$ is negative and $(v_2,v_3)$ is also negative.
Notice that $v_2$ cannot have two siblings as each such sibling is either a leaf or a parent of a leaf, and we will have a connected $8$-edge subgraph with non-positive weight.
Thus, there are two subcases remaining: either $v_2$ has one sibling, or has no sibling. Call the first case $C1$ and call the second case $C2$.

Case C1: Here $v_2$ has one sibling, say $u$. There are two subcases here according to whether $u$ is a leaf, or not. Call the first case $C11$ and call the second case $C12$.

Case C11: If $u$ is a leaf, then $(u,v_3)$ is negative, which forces $(v_3,v_4)$ positive, which forces that $v_3$ has no siblings, which forces $(v_4,v_5)$ positive.
If $v_4$ has no siblings, then $(v_5,v_6)$ must be positive, but then we can remove vertices $v_1,\ldots,v_5,u$ from the tree and use induction on $n-6$ as the total weight hasn't changed
and we remain with a tree. Now, assume $v_4$ has one or more siblings, and say $x$ is such a sibling. Notice that $x$ is in level $t-3$ so to avoid a non-positive $8$-edge subgraph,
we must have all the leaves in the subtree rooted at $x$ in level $t$, and that all edges in this subtree which are non-leaf edges have positive weight. If one of these leaf edges has no
siblings, then this edge, together with the edge incident with it (the edge between level $t-2$ and $t-1$), can be removed and we can use the induction hypothesis with $n-2$.
Otherwise, all leaves of the subtree rooted at $x$ have siblings. The only way this can occur without violating $8$-local positivity is that there are precisely three leaves, say $a,b,c$
(in level $t$) their parent $d$ (in level $t-1$) has no siblings, its parent $f$ (in level $t-2$) has no siblings and the parent of $f$ is $x$.
But then removing $a,b,c,d,f,x$ results in a tree with $n-6$ edges and unchanged weight, hence we can apply induction.

Case C12: Here $u$ is not a leaf. Now $u$ cannot have two or more children (notice that $u$'s children are at level $t$ and hence are leaves), since together with the edges
$(v_1,v_2)$ and $(v_2,v_3)$ we violate $8$-local positivity. Hence, $u$ has exactly one child, say $w$. Now if $(u,v_3)$ is positive then removing $u$ and $w$ we obtain a tree
with $n-2$ edges and unchanged weight, so induction can be applied. Otherwise, if $(u,v_3)$ is negative, then we violate $8$-local positivity.

Case C2: Here $v_2$ has no sibling. There are two subcases here according to whether $v_3$ has siblings, or not. Call the first case $C21$ and call the second case $C22$.

Case C21: let $x$ be a sibling of $v_3$. Then the subtree rooted at $x$ has at most one leaf, since if it had two, we would have two negative edges in this subtree,
and hence a non-positive connected $8$-edge subtree containing these two leaves and $(v_1,v_2)$ and $(v_2,v_3)$. If the subtree rooted at $x$ has exactly one leaf, say $a$, then the path from $a$ to $v_1$
(the path going up from $a$ to $v_4$ and then down from $v_4$ to $v_1$) is a path of length at most $6$
which either contains more than three negative edges, which violates $8$-local positivity, or else contains exactly three negative edges. Now, if this path is of length $6$ or $5$
we can remove $a$ and its parent and use induction with $n-2$ where the total weight is unchanged. If it is of length $4$ then $x$ itself is a leaf. The same reasoning shows that any sibling of $v_3$
in addition to $x$ is a leaf, but this cannot happen since if we had additional siblings then we would violate $8$-local positivity. Thus, $x$ is the only sibling of $v_3$, and $(v_3,v_4)$ is
positive, and hence also $(v_4,v_5)$ is positive. Now, if $v_4$ has a sibling $y$, then in order to not violate $8$-local positivity, all the leaves of the subtree rooted at $y$ have to be in level $t$,
which forcibly enables us to use induction on $n-2$ (if some leaf has no sibling) $n-4$ (if there are two sibling leaves) or $n-6$ (if there are three sibling leaves).
Hence we can now assume that $v_4$ has no sibling, which forces $(v_5,v_6)$ positive. But now, removing vertices $v_1,v_2,v_3,v_4,x,v_5$  we can use induction on $n-6$ since the weight
has not changed.

Case C22: Here $v_3$ has no sibling. We can assume that $(v_3,v_4)$ and $(v_4,v_5)$ are not both positive, since then removing $v_1,v_2,v_3,v_4$ we can use induction on $n-4$ 
since the weight has not changed. We can also assume that $(v_3,v_4)$ and $(v_4,v_5)$ are not both negative, since then the path from $v_1$ to $v_5$ is all negative, violating
$8$-local positivity. So one of $(v_3,v_4)$ $(v_4,v_5)$ is positive and one is negative. But now this means that $v_4$ has no siblings. Indeed, if it had a sibling $x$, then a path from
a leaf in the subtree rooted at $x$ (which can be $x$ itself) to $v_1$ (this path passes through $v_5$ on the way up and then down from $v_5$ to $v_1$) has four negative edges,
which violates $8$-local positivity. So, $v_4$ has no siblings which means that $(v_5,v_6)$ exists and is positive. Now, if $v_5$ has a sibling, say $x$, then the
subtree rooted at $x$ has either four or five levels (since $x$ is in level $t-4$ and the maximum level is $t$, and since all edges between $x$ and its children must be positive and all
edges between these children and their children must be positive). But then it is immediate to check that we can either apply induction on $n-2$, $n-4$, $n-6$ where all removed vertices are in
this subtree, or else violate $8$-local positivity already in this subtree. Hence we can assume that $v_5$ has no sibling so $(v_6,v_7)$ exists and is positive.
But then removing $v_1,\ldots,v_6$  retains a connected tree with $n-6$ edges and the same total weight, so we can apply induction.
\qed

\begin{lemma}\label{l:k-collapse}
$k=7$ and all $k \ge 9$ collapse for $\T$.
\end{lemma}
\Proof
Our constructions depend on the value of $k \pmod 4$. Suppose first that $k=4m+3$. Consider $S_{t,3m+2}$ and assign the last $2m+1$ edges of each branch (the edges closer to the leaves)
the weight $-1$ and the first $m+1$ edges of each branch the weight $1$. This yields a $k$-local positive weighing yet the overall weight is $-mt$. This shows that for every $t \ge 1$ there
is a tree with $n=t(3m+2)$ edges and a $k$-local positive weighing with overall weight $-n \frac{m}{3m+2}$.

If $k=4m+2$ or $k=4m+1$, then consider $S_{t,3m+1}$. Assign the last $2m$ edges of each branch the weight $-1$ and the first $m+1$ edges of each branch the weight $1$.
This yields a $k$-local positive weighing yet the overall weight is $-(m-1)t$. This shows that for every $t \ge 1$ there
is a tree with $n=t(3m+1)$ edges and a $k$-local positive weighing with overall weight $-n \frac{m-1}{3m+1}$.

If $k=4m$, then consider $S_{t,3m}$. Assign the last $2m-1$ edges of each branch the weight $-1$ and the first $m+1$ edges of each branch the weight $1$.
This yields a $k$-local positive weighing yet the overall weight is $-(m-2)t$. This shows that for every $t \ge 1$ there
is a tree with $n=3tm$ edges and a $k$-local positive weighing with overall weight $-n \frac{m-2}{3m}$.
\qed

Given Lemma \ref{l:k-collapse}, it is plausible to suspect that for every $k$, a $k$-local positive weighing of an $n$-edge tree always attains weight at least $-cn$ where $c < 1/3$ 
for all sufficiently large $n$. Indeed, this is true in a much stronger form, and, in fact, in some cases we can prove that the constructions of Lemma \ref{l:k-collapse} are tight.
To this end we need to define the following arithmetic function $g(n,k)$ valid for all positive integers $n,k$.

Suppose $k$ is given. Then define $g(n,k)=-n$ for all $n=1,\ldots,\lfloor (k-1)/2  \rfloor$.
If $k$ is odd then define $g(n,k)=n-k+1$ for all $n= (k+1)/2,\ldots,k$.
If $k$ is even then define $g(n,k)=n-k+2$ for all $n=k/2,\ldots,k$.
Finally, for all $n > k$ define
$$
g(n,k) = \max \{ g(n-1,k)-1 ~,~   \min_{j=\lceil n/3 \rceil}^{\lfloor n/2 \rfloor} g(n-j,k)+g(j,k) \}\;.
$$

\begin{lemma}\label{l:gnk}
Let $k$ and $n$ be positive integers, let $T$ be a tree with at least $k$ edges, and let $w$ be a $k$-local positive weighing of $T$.
Then the weight of every connected $n$-edge subgraph of $T$ is at least $g(n,k)$. In particular, if $T$ has $n$ edges, then $w(T) \ge g(n,k)$.
\end{lemma}
\Proof
Fixing $k$, we prove the lemma by induction on $n$ where the base cases are all $n=1,\ldots,k$. The lemma trivially holds for all $n=1,\ldots,\lfloor (k-1)/2  \rfloor$ since $g(n,k)=-n$ in this case.
Now, if $k$ is odd, any connected subgraph with $n$ edge has at most $(k-1)/2$ negative edges, so the weight of a connected $n$-edge subgraph is at least
$n-2(k-1)/2=n-k+1$. Similarly, if $k$ is even, any connected subgraph with $n$ edge has at most $k/2-1$ negative edges, so the weight of a connected $n$-edge subgraph is at least
$n-2(k/2-1)=n-k+2$. Thus, by the definition of $g(n,k)$, the lemma holds for all $n \le k$.

Now, suppose the lemma holds for all values smaller than $n$ and that $n > k$. Clearly, an $n$-edge connected subgraph has weight at least $g(n-1,k)-1$ by the induction hypothesis for $n-1$.
On the other hand, using again the well-known fact that every tree with more than one edge can be edge-decomposed into two trees with at least $n/3$ edges each, we have
for each $n$-edge subgraph (namely, sub-tree) that there exists $j$ such that $\lceil n/3 \rceil \le j \le \lfloor n/2 \rfloor$ and such that the subgraph can be decomposed into two trees,
on with $j$ edges and the other with $n-j$ edges. But then, by the induction hypothesis used for $j$ and $n-j$, the weight of the subgraph is at least
$g(n-j,k)+g(j,k)$. Thus, the definition of $g(n,k)$ guarantees that the induction step holds.
\qed

In some cases, Lemma \ref{l:gnk}, and the definition of $g(n,k)$ can be used to provide tight bounds for the minimum possible weight of a $k$-local positive weighing
and the asymptotic constant $c_{\T}(k)$.
Table \ref{table:1} shows some initial values of $g(n,k)$ for some representative small $k$.
\begin{table}
\centering
\begin{tabular}{|c||c|c|c|c|c|c|c|c|c|}
\hline
$n \backslash k$ & $7$ & $8$ & $9$ & $10$ &  $11$ & $12$ & $13$ & $14$ & $15$\\
\hline
$1$ & $-1$  & $-1$ & $-1$ & $-1$ & $-1$ & $-1$ & $-1$ & $-1$ & $-1$\\ 
\hline
$2$ & $-2$  & $-2$ & $-2$ & $-2$ & $-2$ & $-2$ & $-2$ & $-2$ & $-2$\\ 
\hline
$3$ & $-3$  & $-3$ & $-3$ & $-3$ & $-3$ & $-3$ & $-3$ & $-3$ & $-3$\\ 
\hline
$4$ & $-2$  & $-2$ & $-4$ & $-4$ & $-4$ & $-4$ & $-4$ & $-4$ & $-4$\\ 
\hline
$5$ & $-1$  & $-1$ & $-3$ & $-3$ & $-5$ & $-5$ & $-5$ & $-5$ & $-5$\\ 
\hline
$6$ & $0$   & $0$  & $-2$ & $-2$ & $-4$ & $-4$ & $-6$ & $-6$ & $-6$\\ 
\hline
$7$ & $1$   & $1$  & $-1$ & $-1$ & $-3$ & $-3$ & $-5$ & $-5$ & $-7$\\ 
\hline
$8$ & $0$   & $2$  & $0$  & $0$  & $-2$ & $-2$ & $-4$ & $-4$ & $-6$\\ 
\hline
$9$ & $-1$  & $1$  & $1$  & $1$  & $-1$ & $-1$ & $-3$ & $-3$ & $-5$\\ 
\hline
$10$ & $-2$ & $0$  & $0$  & $2$  & $0$  & $0$  & $-2$ & $-2$ & $-4$\\ 
\hline
$11$ & $-1$ & $-1$ & $-1$ & $1$  & $1$  & $1$  & $-1$ & $-1$ & $-3$\\ 
\hline
$12$ & $-2$ & $0$  & $-2$ & $0$  & $0$  & $2$  & $0$  & $0$  & $-2$\\ 
\hline
$13$ & $-1$ & $1$  & $-3$ & $-1$ & $-1$ & $1$  & $1$  & $1$  & $-1$\\ 
\hline
$14$ & $-2$ & $0$  & $-2$ & $-2$ & $-2$ & $0$  & $0$  & $2$  & $0$ \\ 
\hline
$15$ & $-3$ & $-1$ & $-3$ & $-1$ & $-3$ & $-1$ & $-1$ & $1$  & $1$ \\ 
\hline
$16$ & $-2$ & $0$  & $-2$ & $0$  & $-4$ & $-2$ & $-2$ & $0$  & $0$ \\ 
\hline
$17$ & $-1$ & $-1$ & $-3$ & $-1$ & $-3$ & $-3$ & $-3$ & $-1$ & $-1$\\ 
\hline
$18$ & $-2$ & $0$  & $-4$ & $-2$ & $-4$ & $-2$ & $-4$ & $-2$ & $-2$\\ 
\hline
$19$ & $-3$ & $1$  & $-3$ & $-1$ & $-3$ & $-1$ & $-5$ & $-3$ & $-3$\\ 
\hline
$20$ & $-4$ & $0$  & $-4$ & $-2$ & $-4$ & $-2$ & $-4$ & $-4$ & $-4$\\ 
\hline
$21$ & $-3$ & $-1$ & $-3$ & $-3$ & $-5$ & $-3$ & $-5$ & $-3$ & $-5$\\ 
\hline
$22$ & $-4$ & $-2$ & $-2$ & $-2$ & $-4$ & $-2$ & $-4$ & $-2$ & $-6$\\ 
\hline
$23$ & $-3$ & $-1$ & $-3$ & $-1$ & $-5$ & $-3$ & $-5$ & $-3$ & $-5$\\ 
\hline
$24$ & $-4$ & $0$  & $-4$ & $0$  & $-6$ & $-4$ & $-6$ & $-4$ & $-6$\\ 
\hline
$25$ & $-5$ & $-1$ & $-5$ & $-1$ & $-5$ & $-3$ & $-5$ & $-3$ & $-5$\\ 
\hline
$26$ & $-4$ & $-2$ & $-6$ & $-2$ & $-4$ & $-4$ & $-6$ & $-4$ & $-6$\\ 
\hline
$27$ & $-5$ & $-1$ & $-5$ & $-3$ & $-5$ & $-3$ & $-7$ & $-5$ & $-7$\\ 
\hline
$28$ & $-4$ & $-2$ & $-6$ & $-4$ & $-4$ & $-2$ & $-6$ & $-4$ & $-6$\\ 
\hline
$29$ & $-5$ & $-1$ & $-5$ & $-3$ & $-5$ & $-1$ & $-7$ & $-5$ & $-7$\\ 
\hline
$30$ & $-6$ & $-2$ & $-6$ & $-2$ & $-6$ & $-2$ & $-6$ & $-6$ & $-8$\\ 
\hline
$31$ & $-5$ & $-1$ & $-7$ & $-3$ & $-7$ & $-3$ & $-5$ & $-5$ & $-7$\\ 
\hline
$32$ & $-6$ & $-2$ & $-6$ & $-4$ & $-8$ & $-4$ & $-6$ & $-4$ & $-8$\\ 
\hline
$33$ & $-5$ & $-3$ & $-7$ & $-3$ & $-7$ & $-5$ & $-5$ & $-3$ & $-9$\\
\hline
$34$ & $-6$ & $-2$ & $-6$ & $-4$ & $-8$ & $-6$ & $-6$ & $-2$ & $-8$\\ 
\hline
$35$ & $-7$ & $-1$ & $-7$ & $-5$ & $-7$ & $-5$ & $-7$ & $-3$ & $-7$\\
\hline
\end{tabular}
\caption{Some values of $g(n,k)$.}
\label{table:1}
\end{table}

\begin{lemma}\label{l:bounds}
Suppose for some integer $r$ we have that $3r > k$ and $g(n,k)/n \ge c$ for all $n=r,\ldots,3r$.
Then for any $k$-local positive weighing of a tree with $n \ge k$ edges it holds that $w(T) \ge cn$.
In particular, if $k$ collapses for $\T$, then $c_T(k) \ge c$.
Consequently, it holds that
\begin{enumerate}
\item
If $k = 4m+3$, then $c_\T(k) = -\frac{m}{3m+2}$\;.
\item
If $k = 4m+2$, then $c_\T(k) = -\frac{m-1}{3m+1}$\;.
\item
If $k = 4m+1$, then $c_\T(k) \ge -\frac{2m-1}{6m+1}$\;.
\item
If $k = 4m$, then $c_\T(k) \ge -\frac{2m-3}{6m-1}$\;.
\end{enumerate}
\end{lemma}
\Proof
By Lemma \ref{l:gnk}, it suffices to prove that for all $n \ge r$ we have that $g(n,k) \ge cn$.
Indeed by our assumption, this holds for all $n=r,\ldots,3r$. Assume now that $n > 3r$ and that $g(j,k) \ge cj$ holds for all $j$
with $r \le j < n$. We prove that also $g(n,k) \ge cn$. Indeed, since $n > 3r > k$, we have by the definition of $j$ that
for some $j$  with $\lceil n/3 \rceil \le j \le \lfloor n/2 \rfloor$ it holds that $g(n,k) \ge g(j,k)+g(n-j,k)$.
But both $j$ and $n-j$ are less than $n$ and larger than $r$. Therefore,  $g(n,k) \ge cj+c(n-j)=cn$.

Having proved the first part of the lemma, we demonstrate its use in some particular cases.
For $k=7$ we see from Table \ref{table:1} that for $r=5$, $g(n,7)/n \ge -1/5$ for all $n=5,\ldots,15$.
Hence $C_\T(7) \ge -1/5$.  Equality holds by the construction in Lemma \ref{l:k-collapse}.
For $k=10$ we see from Table \ref{table:1} that for $r=7$, $g(n,10)/n \ge -1/7$ for all $n=7,\ldots,21$.
Hence $C_\T(10) \ge -1/7$.  Equality holds by the construction in Lemma \ref{l:k-collapse}.

We next prove that when $k=4m+2$ , then $c_\T(k) \ge -\frac{m-1}{3m+1}$.
Note that if we do that, then equality holds by the construction in Lemma \ref{l:k-collapse}.
By the first part of the lemma, it suffices to show that for, say, $r=3m+1$ we have that $g(n,k)/n \ge -\frac{m-1}{3m+1}$
for all $n=r,\ldots,3r$.
Indeed this holds true initially for $g(r,k)$ since in this case $k$ is even and $r$ falls
in the range where the function's definition is $g(r,k)=r-k+2=-m+1$. So, $g(r,k)/r = -\frac{m-1}{3m+1}$, as required.
Now notice that $g(n,k)$ increases each time by $1$ for $n=r+1,\ldots,k$ where it reaches its maximum $g(k,k)=2$ and then starts to decrease by $1$
until $n=2r$ where it equals $g(2r,k)=-2m+2$. Hence the minimum of $g(n,k)/n$ in the range $n=r+1,\ldots,2r$ is attained for $n=2r$ where it is
again $-\frac{m-1}{3m+1}$, as required.
Now we need to consider the range $n=2r, \ldots 3r$. By the definition of $g(n,k)$ it suffices to prove that for each
$j$ such that $\lceil n/3 \rceil \le j \le \lfloor n/2 \rfloor$ it holds that
$$
\frac{g(j,k) + g(n-j,k)}{n} \ge -\frac{m-1}{3m+1}\;.
$$
We consider three cases according to the value of $j$.

First case: $j \le k$ and $j \ge n-k$.
Now, since $g(r,k)=-m+1$, then $g(j,k) \ge -m+1+j-r$ since the absolute difference between any two consecutive values of $g(.,k)$ is $1$.
On the other hand, since we assume $k \ge n-j$, we have that $g(n-j,k) \ge g(k,k) - k + n - j = 2-k-j+n$.
Hence
$$
g(j,k) + g(n-j,k) \ge -m+1+j-r + 2-k-j+n = -m+3-r-k+n = -8m+n\;.
$$
Now, $(-8m+n)/n$ is minimized when $n$ is minimized, namely when $n=2r=6m+2$, in which case the ratio is $-\frac{m-1}{3m+1}$, as required.

Second case: $j \le k$ and $j \le n-k$.
Again we use $g(j,k) \ge -m+1+j-r$. But now, since $k \le n-j$, we have that $g(n-j,k) \ge g(k,k) - (n-j) + k = 2+k+j-n$.
Hence
$$
g(j,k) + g(n-j,k) \ge -m+1+j-r + 2+k+j-n = -m+3-r+k-n+2j = 4-n+2j \ge 4-n/3\;.
$$
Now, $(4-n/3)/n$ is minimized when $n$ is maximized, namely when $n=3r=9m+3$, in which case the ratio is at least $-\frac{m-1}{3m+1}$, as required.

Third case: $j \ge k$. 
We have
$$
g(j,k) + g(n-j,k) \ge g(k,k)-(j-k) + g(k,k) - ((n-j)-k) = 4-n+2k\;.
$$
Now, $(4-n+2k)/n$ is minimized when $n$ is maximized, namely when $n=3r$. But even in this case we have
$$
\frac{4-3r+2k}{3r} = \frac{5-m}{3+9m} > -\frac{m-1}{3m+1}\;.
$$
The other three cases of $k=4m+3$, $k=4m+1$ and $k=4m$ are proved analogously. Notice that for the case
$k=4m+3$ we obtain equality by the construction in Lemma \ref{l:k-collapse}.
\qed

\vspace{5pt}
\noindent
{\bf Proof of Theorem \ref{t:1}} now follows from Theorem \ref{t:trees-graphs} and Lemmas \ref{l:1234}, \ref{l:k56}, \ref{l:k8}, \ref{l:k-collapse}, \ref{l:bounds}.
\qed

\section{Bounded degree graphs}

Recall that for $\Delta \ge 2$, $\G_\Delta$ denotes the subset of a graph family $\G$ consisting of all elements with maximum degree at most $\Delta$.
Our aim in this section is to classify $k$ for all $\G_\Delta \supseteq \T_\Delta$. This classification differs somewhat from
the general case given in Theorem \ref{t:1}.
First observe that every $k \ge 1$ is forcing for $\G_2=\T_2$ (namely, for the family of all paths and all cycles).
This follows trivially from the fact that a path or cycle of length $n > k$ which is $k$-positively weighted has total weight larger than $n/k-k$.
So from this point onwards, we are only concerned with the case $\Delta \ge 3$.

\begin{lemma}\label{l:delta-k1234}
For all $k=1,2,3,4$ and for all $\Delta \ge 3$ we have that $k$ is forcing for $G_\Delta$.
\end{lemma}
\Proof
This immediately follows from the fact that $k=1,2,4$ are forcing already for $\G$ and from Lemma \ref{l:1234} 
stating that for every $\Delta$, $k=3$ is forcing for $\T_\Delta$. Notice that for $k=3$, forcing for $\T_\Delta$ implies forcing for $\G_\Delta$ by 
Theorem \ref{t:trees-graphs}. \qed

\begin{lemma}\label{l:delta-k6}
$k=6$ is forcing for $G_\Delta$ for all $\Delta \ge 3$.
\end{lemma}
\Proof
First observe that by Theorem \ref{t:trees-graphs} we have that if $6$ is forcing for $\T_\Delta$, then it is forcing for $\G_\Delta$.
So, we have to prove that $6$ is forcing for $\T_\Delta$. We will prove that for every tree in $T_\Delta$ with more than $4\Delta$ edges,
a $6$-local positive weighing implies $w(T) > 0$. Observe that the requirement  of having more than $4\Delta$ edges is tight,
since $S_{\Delta,4} \in T_\Delta$ is a tree with $4\Delta$ edges and the weighing that assigns $-1$ to the leaf edges and to the edges adjacent to them
and $1$ to all other edges is $6$-local positive, while its total weight is $0$.

We first notice that by the proof in Lemma \ref{l:k56}, every $6$-local positive weighing of any tree with at least $6$ edges has nonnegative weight.
So, assume now that $T \in T_\Delta$ has more than $4\Delta$ edges. We use the same rooted tree notations as in Lemma \ref{l:k56}
and recall that we may assume that all leaf edges are negative, as otherwise we can delete a leaf edge and the remaining tree
(which still has more than $6$ edges) has nonnegative weight, hence $T$ has positive weight.
As in Lemma \ref{l:k56}, we have two cases according to whether $v_1$ has a unique sibling, or has no siblings.
We prove for the case that $v_1$ has a unique sibling $v_0$. The proof of the other case is analogous.

So, $(v_1,v_2)$ and $(v_0,v_2)$ are both leaf edges and are both negative. Hence $(v_2,v_3)$ is positive and as in the proof of Lemma \ref{l:k56},
$v_2$ has no sibling and $(v_3,v_4)$ is positive. However, unlike Lemma \ref{l:k56}, we would not gain much by deleting vertices
$v_0,v_1,v_2,v_3$ since we remain with a tree on $n-4$ edges and unchanged total weight, but this would only tell us that $w(T) \ge 0$,
which is not enough. Instead, we argue as follows. Let $u$ be any child of $v_4$ (so, in particular, $v_3$ is such a child, but there may be others).
Then $(u,v_4)$ is positive so $u$ is not a leaf, and also no child of $u$ is a leaf (as otherwise, together with the negative edges $(v_1,v_2)$ and $(v_0,v_2)$
we would have a connected non-positive $6$-edge subgraph). So, either the sub-tree rooted at $u$ has (i) two edges, one positive and one negative,
or (ii) has three edges, one positive and two negative, or (iii) has four edges, two positive and two negative. In cases (i) and (iii) we can remove
the subtree rooted at $u$ (which also removes $u$ and the edge $(u,v_4)$) and we have removed a positive weighted tree on at most $5$ edges,
and the remaining tree has at least $4\Delta-5 \ge 7$ edges, hence has nonnegative total weight, and we are done.

We remain with case (ii) and we can assume that every child $u$ of $v_4$ satisfies (ii) (notice that $v_3$ is indeed such a child).
So, the whole subtree rooted at $v_4$ has $4x$ edges and total weight $0$, where $x$ is the number of children of $v_4$.
Since $T\ \in \T_\Delta$ we have $x \le \Delta$. We claim that $x < \Delta$. Indeed, otherwise $v_4$ is the root and the whole tree has
$4\Delta$ edges, while we assume it has more. So, $x \le \Delta-1$ and once removing all edges of the subtree rooted at $v_4$ (but not $v_4$
itself), we remain with a tree on at least $n-4(\Delta-1) \ge 5$ edges. But in the new tree, $v_4$ is a leaf and the leaf edge incident to
it (namely $(v_4,v_5)$) has positive weight as otherwise, together with the negative edges $(v_1,v_2)$ and $(v_0,v_2)$
we would have a connected non-positive $6$-edge subgraph.  But then we can remove this leaf, and remain with a nonnegative weighted tree,
as it has at least $4$ edges (if it has at least $6$ edges, it has nonnegative weight by assumption, and if it has $4$ or $5$ edges, it has at most
two negative edges, hence is also nonnegative). Hence $T$ has positive weight. \qed

\begin{lemma}\label{l:delta-k58}
For all $\Delta \ge 3$, each of $k=5,8$ is weakly forcing for $\G_\Delta$ with $f_{\G_\Delta}(5)=-1$, $f_{\G_\Delta}(8)=-1$.
\end{lemma}
\Proof
First observe that by Theorem \ref{t:trees-graphs} it suffices to prove that $5$ and $8$ are weakly forcing for $\T_\Delta$ and 
$f_{\T_\Delta}(5)=-1$, $f_{\T_\Delta}(8)=-1$. In fact, it suffices to prove that $f_{\T_3}(5) \le -1$ and $f_{\T_3}(8) \le -1$
since $\T_3 \subseteq \T_\Delta \subset \T$ and we already have the other direction from Lemmas \ref{l:k56} and \ref{l:k8}.

Let $C_{t,r}$ denote the so-called {\em comb} which is obtained by taking a path on $t$ vertices (this is the central path) and attaching disjoint paths of
length $r$ to each vertex of the central path. Notice that $C_{t,r}$ has $rt+t-1$ edges and that $C_{t,r} \in \T_3$.
Now specifically consider $C_{t,3}$ and assign weight $-1$ to all leaf edges and the edges adjacent to them and weight $1$ to the other edges.
This constitutes a $5$-local positive weighing with total weight $-1$, proving $f_{\T_3}(5) \le -1$.
Now take $C_{t,5}$ and assign weight $-1$ to all leaf edges, to the edges adjacent to them, and to the edges adjacent to them and weight $1$ to the other edges.
This constitutes an $8$-local positive weighing with total weight $-1$, proving $f_{\T_3}(8) \le -1$.
\qed

\begin{lemma}\label{l:delta-k-collapse}
For all $\Delta \ge 3$ we have that $k=7$ and all $k \ge 9$ collapse for $\G_\Delta$.
\end{lemma}
\Proof
We will prove collapse already for $\T_3$, hence also for $\G_\Delta$ which contains it.
Our constructions depend on the value of $k \pmod 4$.
For the case $k = 4m$ take $C_{t,3m-1}$ and assign weight $-1$ to the $2m-1$ edges closest to each leaf,
and weight $1$ to the other edges. This constitutes a $k$-local positive weighing with total weight $-(m-2)t-1$.
This shows that for every $t \ge 1$ there
is a tree in $\T_3$ with $n=3mt-1$ edges and a $k$-local positive weighing with overall weight $-(m-2)\frac{n+1}{3m}-1$.

For the case $k=4m+1$, take $C_{t,3m}$ and assign weight $-1$ to the $2m$ edges closest to each leaf,
and weight $1$ to the other edges. This constitutes a $k$-local positive weighing with total weight $-(m-1)t-1$.
This shows that for every $t \ge 1$ there
is a tree in $\T_3$ with $n=(3m+1)t-1$ edges and a $k$-local positive weighing with overall weight $-(m-1)\frac{n+1}{3m+1}-1$.

For the cases $k=4m+2$ and $k=4m+3$ we need an alternating comb, which we denote by $A_{t,r}$ and defined as follows.
Take a path on $t$ vertices (this is the central path) and attach disjoint paths of
length $r$ to each vertex of the central path which is located in an odd place on the central path and
disjoint paths of
length $r+1$ to each vertex of the central path which is located in an even place on the central path.

For the case $k=4m+2$, take $A_{2q+1,3m}$. The $2m$ edges closest to each leaf get weight $-1$ and the others get weight $1$.
This constitutes a $k$-local positive weighing with total weight $-q(2m-3)-m$.
This shows that for every $q$, there
is a tree in $\T_3$ with $n=(6m+3)q+3m$ edges and a $k$-local positive weighing with overall weight $-(2m-3)\frac{n-3m}{6m+3}-m$.

For the case $k=4m+3$, take $A_{2q+1,3m+1}$. The $2m+1$ edges closest to each leaf get weight $-1$ and the others get weight $1$.
This constitutes a $k$-local positive weighing with total weight $-q(2m-1)-m-1$.
This shows that for every $q$, there
is a tree in $\T_3$ with $n=(6m+5)q+3m+1$ edges and a $k$-local positive weighing with overall weight $-(2m-1)\frac{n-3m-1}{6m+5}-m-1$.
\qed

\vspace{5pt}
\noindent
{\bf Proof of Theorem \ref{t:2}} now follows from Lemmas \ref{l:delta-k1234}, \ref{l:delta-k6}, \ref{l:delta-k58},
\ref{l:delta-k-collapse}, and the lower bounds for $c_{\G}(k)$, given in Theorem \ref{t:1}.
\qed

\section{Graphs with bounded average degree}

\begin{lemma}\label{l:avg-1}
For every $\alpha < 8$, all sufficiently large $k$ collapse for $\C_{\overline{\alpha}}$.
\end{lemma}
\Proof
We will assume that $k \equiv 3 \pmod 4$. The other cases are proved analogously.
Consider a path $P$ of length $(k+1)/4$, so it has $(k+5)/4$ vertices. Denote one endpoint of the path by $x$ and the other endpoint by $y$.
Now let $q$ be the least integer such there exists a graph $Q$ on $q$ vertices which contains $3(k-3)/4$ edges and
observe that $q = O(\sqrt{ k})$. Attach a copy of $Q$ to $P$ by identifying one vertex of $Q$ with $y$.
Call the new graph $R$. Now observe that $R$ has $(k+1)/4+q$ vertices and $(k+1)/4+3(k-3)/4=k-2$ edges. 
Finally, take $t$ copies of $R$ and identify all their $x$ vertices. Denote this graph by $H_{k,t}$.

Now, $H_{k,t}$ has $t((k-3)/4+q)+1$ vertices and $t(k-2)$ edges. Its average degree is therefore
$$
\frac{2t(k-2)}{t((k-3)/4+q)+1}\;.
$$
Hence, for every $\alpha < 8$ we have that for all $k$ sufficiently large, $H_{k,t} \in \C_{\overline{\alpha}}$ for every $t$.

It remains to show that there is a $k$-local positive weighing of $H_{k,t}$ whose total weight is $-t$,
thereby proving that all $k$ sufficiently large collapse for $\C_{\overline{\alpha}}$.

To each edge that belongs to $P$ assign weight $1$. To exactly $(k-1)/2$ edges of $Q$ assign weight $-1$ and to the remaining
$3(k-3)/4-(k-1)/2=(k-7)/4$ edges off $Q$ assign weight $1$. It is immediate to verify that
this is a $k$-local positive weighing of $H_{k,t}$ since the distance between two $Q$ blocks must pass through two $P$ blocks completely.
Yet the total weight of each of the $t$ copies of $R$ is
$$
-\frac{k-1}{2} + \frac{k+1}{4} + \frac{k-7}{4} = -1\;.
$$
Hence the weight of $H_{k,t}$ is $-t$,
\qed

\begin{lemma}\label{l:avg-2}
For every $\alpha \ge 8$, all $k$ are forcing for $\C_{\overline{\alpha}}$.
\end{lemma}
\Proof
Let $G \in C_{\overline{8}}$ be a given $n$-vertex graph with a $k$-local positive weighing and notice that $m=|E(G)| \ge 4n$.
Our goal is to prove that more than $m/2$ edges have positive wight.
We will assume that $n$ is sufficiently large as a function of $k$, as the definition of forcing allows for finitely many exceptions.

As long as there is a positive edge which is on a cycle, we remove it from $G$. We halt this process once there are no more such edges.
Suppose we have removed $x$ edges. Then, we remain with a connected graph $G'$ on $m-x$ edges which still satisfies $k$-local positivity,
and $w(G)=w(G')+x$. Now, by Theorem \ref{t:1}, we know that $w(G') > -(m-x)/3$. hence, $w(G) > -m/3+4x/3$.
It therefore suffices to prove that $x \ge m/4$.

Since $w(G') > -(m-x)/3$, the number of positive edges of $G'$ is larger than $(m-x)/3$.
Now, since all positive edges of $G'$ are bridges in $G'$, then the number of bridges is larger than $(m-x)/3$. 
Notice that if $m \ge 4n$, then we immediately obtain that $x \ge m/4$ as otherwise $(m-x)/3 \ge n$ which is impossible
as the number of bridges in an $n$-vertex graph cannot exceed $n-1$. 
\qed

\vspace{5pt}
\noindent
{\bf Proof of Theorem \ref{t:3}} now follows from Lemmas \ref{l:avg-1} and \ref{l:avg-2}.
\qed

\section{Concluding remarks}

Theorem \ref{t:1} classifies all $k$ for every connected family of graphs $\G$ that contains all trees.
In the collapsing case, the precise value of $c_{\G}(k)$ is determined for two moduli of $4$, while for two
other moduli, only upper and lower bounds are provided (although they do approach each other as $k$ grows).
It would be of some interest to determine $c_{\G}(k)$ in these remaining cases as well,
and the same holds for the bounds for $c_{\G_\Delta}(k)$ in Theorem \ref{t:2}.
At least for small $k$, it is possible to show that in the case where an exact bound is not determined, the lower bounds are not tight.
Indeed, say, for $k=9$ it is not difficult to show that a $9$-local positive tree with $26$ edges has total weight at least $-4$.
Hence we can define a variant of $g(n,9)$ (the function from Table \ref{table:1}) where we set $g(26,9)=-4$ (instead of the present value of $-6$)
and continue with the same recursive definition. It is simple to verify that the revised $g(n,9)/n$ never falls below $-7/31$ for all $n \ge 14$,
which, by Lemma \ref{l:gnk} implies that $c_{\G}(9) \ge -7/31$.

Another variant of $g(n,k)$, where instead of defining $g(k,k)=2$ for even $k$ we define $g(k,k)=0$ (instead of the present value of $+2$)
and continue with the same recursive definition, can be used to prove that if the total weight of a sufficiently large $m$-edge tree falls below $-m/3$
then it must have a {\em negative} connected subgraph with $k$ edges. This is tight when $k = 0 \pmod 4$ as a construction
similar to the one used in \ref{l:k-collapse} shows that we can have total weight $-m/3$ while every connected $k$-edge subgraph is nonnegative.

Finally, it may be interesting to classify $k$ for other, less general (but perhaps more structured) families of infinite graphs
or for other types of combinatorial objects where a notion of connectivity is meaningful.

\bibliographystyle{plain}

\bibliography{paper}

\end{document}